\documentclass{ifacconf}

\usepackage{graphicx}      

\usepackage{amssymb,amsmath}
\usepackage{epsfig} 
\usepackage{color}
\usepackage{caption,subcaption}
\usepackage{nicefrac}

\makeatletter
\let\old@ssect\@ssect 
\makeatother

\usepackage{natbib}        
\usepackage[colorlinks=true,breaklinks=true,bookmarks=true,urlcolor=blue,citecolor=blue,linkcolor=blue,bookmarksopen=false,draft=false]{hyperref}

\makeatletter
\def\@ssect#1#2#3#4#5#6{%
	\NR@gettitle{#6}
	\old@ssect{#1}{#2}{#3}{#4}{#5}{#6}
}
\makeatother

\newtheorem{theorem}{\bf Theorem}
\newtheorem{proposition}{\bf Proposition}
\newtheorem{definition}{\bf Definition}
\newtheorem{lemma}{\bf Lemma}
\newtheorem{corollary}{\bf Corollary}

\newtheorem{remark}{\bf Remark}

\newcommand{\dfb}{\stackrel{\Delta}{=}}

\def\E{\mathbb{E}}

\def\vec{\mathop{\rm vec}\nolimits}

\def\calV{\mathcal{V}}

\DeclareMathOperator*{\argmax}{argmax} 

\def\qed{\hfill$\blacksquare$}

\def\one{\mathbf{1}}

\newcommand{\ups}{\upsilon}
\newcommand{\ov}{\overline}
\newcommand{\mc}{\mathcal}

\newcommand{\R}{\mathbb{R}}
\newcommand{\ba}{\begin{array}}
\newcommand{\ea}{\end{array}}
\newcommand{\beq}{\begin{equation}}
\newcommand{\eeq}{\end{equation}}
\newcommand{\beqn}{\begin{equation*}}
\newcommand{\eeqn}{\end{equation*}}
\newcommand{\be}{\begin{equation}}
\newcommand{\ee}{\end{equation}}
\DeclareMathOperator*{\argmin}{argmin}

\begin{document}
\begin{frontmatter}



\title{Reaching Optimal Distributed Estimation Through Myopic Self-Confidence Adaptation
}
\thanks[footnoteinfo]{This work was partly supported by the Italian Ministry for University and Research through grants ``Dipartimenti
di Eccellenza 2018–2022'' [CUP: E11G18000350001] and Project PRIN 2017 ``Advanced Network Control of Future Smart Grids''
(http://vectors.dieti.unina.it), and by the Compagnia di San Paolo.}
\author[DISMA]{Giacomo Como}
\author[DISMA]{Fabio Fagnani}
\author[DET]{Anton V. Proskurnikov}
\address[DISMA]{Department of Mathematical Science ``G.L. Lagrange'',  Politecnico di Torino, 10129 Torino, Italy. \\
	E-mail: \texttt{\{giacomo.como, fabio.fagnani\}@polito.it.}}
\address[DET]{Department of Electronics and Telecommunications,  Politecnico di Torino, 10129 Torino, Italy. E-mail: \texttt{anton.p.1982@ieee.org}}

\begin{abstract}
Consider discrete-time linear distributed averaging dynamics, whereby a finite number of agents in a network start with uncorrelated and unbiased noisy measurements of a common state of the world modeled as a scalar parameter, and iteratively update their estimates following a non-Bayesian learning rule. Specifically, let every agent update her estimate to a convex combination of her own current estimate and those of  her neighbors in the network
(this procedure is also known as the French-DeGroot model, or the consensus algorithm).
As a result of this iterative averaging process, each agent obtains an asymptotic estimate of the state of the world, and the variance of this individual estimate depends on the matrix of weights the agents assign to themselves and to the others. We study a game-theoretic multi-objective optimization problem whereby every agent seeks to choose her self-confidence value in the convex combination in such a way to minimize the variance of her asymptotic estimate of the state of the world. Assuming that the relative influence weights assigned by the agents to their neighbors in the network remain fixed and form an irreducible relative influence matrix, we characterize the Pareto frontier of the problem, as well as the set of Nash equilibria in the resulting game.
\end{abstract}

\begin{keyword}
Games on graphs, centrality measures, opinion dynamics.\\
\emph{2010 MSC: 05C57,91A43,91D30}
\end{keyword}

\end{frontmatter}

\section{Introduction}

The wisdom of crowds~\citep{Surowiecki:2004} is a well studied phenomenon. Its central statement is that cooperative decision making generally outperforms individual one. A typical mathematical formalization is the context of a group of agents collecting independent noisy measures of the same quantity. Basic inferential arguments imply the existence of linear aggregative functions of such data that yield a reduced noise variance over the individual ones.

When the aggregation mechanism is obtained through an influence system, however, the final outcome might in principle lead to suboptimal estimations and, in certain extreme situations, the wisdom of crowds effect might vanish. This happens for instance when all agents have equally good estimations, but the influence system has a bias towards some agents and as a consequence it does not lead to a uniform aggregation. More generally, this can occur when estimations are of different quality and the influence system does not properly incorporate this information in the aggregation mechanism. See, e.g., \cite{GolubJackson:2010}, \cite{Bullo.Fagnani.Franci:2020}, and references therein.

In this paper we consider one of the simplest models for the influence system shaping the opinions of a set of agents interacting in a social network, first introduced in \cite{French:1956},~\cite{DeGroot} and~\cite{Lehrer:1976}. See also~\cite{ProTempo:2017-1} and~\cite{BulloBook-Online} for updated accounts of the large body of literature on its analysis. In spite of its simplicity, the French-DeGroot model is more than a purely theoretic construct; its predictive power, for instance, has been recently demonstrated by experiments with large-scale social networks~\citep{Kozitsin:2020}.
The agents' initial states consist of a family of uncorrelated unbiased noisy measures of a common state of the world, modeled as an unknown scalar value. We consider the context where interpersonal influence weights have already been assigned and the agents are only free to modify the diagonal part of the influence matrix (modifying their self-confidence levels) and rescaling the rest. If the agents could directly cooperate, a prescribed assignments of self-confidence weights could be chosen to obtain the optimal final estimation for all of them.

In this paper, we consider a game theoretic scenario where single agents myopically modify self-confidence to maximize the precision of their final estimation. The peculiarity of this game is that the utility functions present a discontinuity when any of the self-confidence parameters is assigned value $1$, as this disconnects the graph and dramatically changes the nature of the final outcome.

Our main contribution consists in a detailed analysis of the pure strategy Nash equilibria of this game. Our main results are stated in Theorems~\ref{theo:main} and~\ref{theo:other-Nash}. Essentially, our analysis  shows that the considered game always admits pure strategy Nash equilibria that correspond to self-confidence choices that do not disconnect the graph. All these equilibria are strict and equivalent in the sense that, under them, all agents  obtain the same best possible estimation, which coincides with the  one they would have obtained by an optimal direct cooperation.
For generic choices of the variance values, no other Nash equilibrium exists. In contrast, when there exist subsets of at least two agents with the same variance other Nash equilibria may appear.  Such Nash equilibria always arise due to disconnection of a number of agents with the same variance from the rest of the agents and are never strict. Moreover, our analysis shows two further insights. First, finite sequences of best response actions exist from Nash equilibria of this second type that lead to a fully connected configuration. Second, in any fully connected configuration, best response actions will never disconnect the graph.

We conclude with a brief outline of the paper. In Section \ref{sec:start} we present the model and an overview of the main result. Section \ref{sec:main} is the main technical section, which develops a fundamental analysis of the best response sets of this game; the proofs of main results are given in Section~\ref{sec.proofs}.
In Section~\ref{sec:dynam}, we present a dynamical system (being a modification of the gradient-type best response dynamics) and some numerical simulations illustrating its behavior.

\section{Problem setup}\label{sec:start}

\subsection{Notational convections and graph-theoretic notions}\label{sec:notation}
Unless otherwise stated, all vectors are considered as columns. We use $\one$ to denote the column vector of ones (its size is being clear from the context).
For two vectors $x$ and $y$ in $\mathbb{R}^n$, the inequalities $\leq,<,\geq,>$ are meant to hold true entry-wise.
For a vector $z$ in $\mathbb{R}^n$, the symbol $[z]$ denotes the diagonal $n\times n$ matrix with diagonal entries $z_1,\ldots,z_n$.
Throughout the paper, the transpose of a matrix $M$ is denoted by $M'$. As usual, a \emph{row-stochastic} matrix is a nonnegtive square matrix $M$ in $\R_+^{n\times n}$ whose row sums are all equal to $1$, i.e., such that $M\one=\one$.

A (finite, directed) \emph{graph} is the pair $\mc G=(\mc V,\mc E)$, of a finite set of nodes $\mc V$ and a set of (directed) links $\mc E\subseteq\mc V\times\mc V$.
A graph $\mc G=(\mc V,\mc E)$ is called \emph{undirected} if $(i,j)\in\mc E$ if and only if $(j,i)\in\mc E$.
The \emph{in-degree} and \emph{out-degree} of a node $i$ in a graph $\mc G$ are defined as
$d_i^-=|\{j\in\mc V:\,(j,i)\in\mc E\}|$  and  $d_i^+=|\{j\in\mc V:\,(j,i)\in\mc E\}|$, respectively. For $k\ge0$, graph $\mc G=(\mc V,\mc E)$ is $k$-regular if $d_i^-=d_i^+=k$ for every $i$ in $\mc V$.

A length-$l$ \emph{walk} from a node $i$ to a node $j$ in a graph $\mc G=(\mc V,\mc E)$ is an $(l+1)$-tuple of nodes $(\gamma_0,\gamma_1,\ldots,\gamma_l)$ such that $\gamma_0=1$, $\gamma_l=j$, and $(\gamma_{k-1},\gamma_{k})\in\mc E$ for  $1\le k\le l$. A walk $(\gamma_0,\gamma_1,\ldots,\gamma_l)$ is referred to as a \emph{path} if $\gamma_h\ne\gamma_k$ for every $0\le h<k\le l$, except for possibly $\gamma_0=\gamma_l$, in which case the path is referred to as a \emph{cycle}.

For a graph $\mc G=(\mc V,\mc E)$ and a subset $\mc S\subseteq\mc V$, we let $\mc G[\mc S]=(\mc S, \mc E[\mc S])$ be the graph with node set $\mc S$ such that, for every $i\ne j$ in $\mc S$, $(i,j)\in\mc E[\mc S]$ if and only if there exists a path from $i$ to $j$ in $\mc G$ that does not pass through any intermediate node $k$ in $\mc S$ (in particular, if $(i,j)\in\mc E$ and $i,j\in\mc S$, then $(i,j)\in\mc E[\mc S]$).

A graph $\mc G=(\mc V,\mc E)$ is referred to as:  \emph{strongly connected} if for every two nodes $i\ne j$ in $\mc V$ there exists a path from $i$ to $j$; \emph{aperiodic} if the greatest common divisor of the lengths of all its cycles is equal to $1$.
A \emph{directed ring} is a $1$-regular, strongly connected graph $\mc G=(\mc V,\mc E)$.

Finally, to a nonnegative square matrix $M$ in $\R_+^{n\times n}$ we can associate the graph $\mc G_M=(\mc V,\mc E)$ with node set $\mc V=\{1,\ldots,n\}$ and link set $\mc E$ such that $(i,j)\in\mc E$ if and only if $M_{ij}>0$. A nonnegative square matrix $M$ in $\R_+^{n\times n}$ is \emph{irreducible} if its associated graph $\mc G_M$ is strongly connected, and \emph{aperiodic} if $\mc G_M$ is aperiodic.

\subsection{Opinion consensus dynamics and social power}

Consider a finite set of agents $\calV=\{1,\ldots,n\}$, each characterized by a real scalar \emph{opinion} $x_i$ and a self-confidence value $z_i$ in $[0,1]$. Social influence between the agents is captured by a row-stochastic matrix $P$ in $\R^{n\times n}$, to be referred to as the \emph{influence matrix}. Without loss of generality we may assume that $P$ has zero diagonal. The influence matrix $P$ defines the social network, whose topology is identified by the graph $\mc G_P$.

Upon stacking all the agents' opinions and self-confidence values respetively in vectors $x$ (the opinion profile) in $\R^n$ and $z$ (the self-confidence profile) in $\mc Z=[0,1]^n$, we may compactly write the French-DeGroot opinion dynamics model as the discrete-time system
\beq\label{eq.degroot}
x(t+1)=((I-[z])P+[z])x(t)\,,\qquad t=0,1,\ldots\,.
\eeq
The above recursion entails every agent to update her opinion to a convex average
$$x_i(t+1)=z_ix_i(t)+(1-z_i)\sum_{j=1}^nP_{ij}x_j(t)$$
of her own current opinion and her neighbors' ones. Specifically, agent $i$'s self-confidence value $z_i$ represents the weight that agent $i$ puts on her own current opinion $x_i(t)$, while each of her neighbors' current opinions $x_j(t)$ receives a weight equal to the product of the  of agent $i$'s complementary self-confidence $(1-z_i)$ and the relative influence weight $P_{ij}$. We shall refer to an agent $i$ as stubborn if $z_i=1$; a stubborn agent keeps her opinion constant in time regardless of the opinion profile.
We shall denote by
\be
\mc S(z)=\{i=1,\ldots,n:\,z_i=1\}\,,\ee
the set of all stubborn agents.

Throughout, we shall restrict our analysis to social networks whose influence matrix $P$ is irreducible.
It then follows from classical Perron-Frobenius theory~\cite[Ch.2, Th.~1.3]{BermanPlemmons_Book} that in this case there exists a unique invariant probability distribution
\be\label{pi}\pi=P'\pi\,,\qquad \one'\pi=1\,,\ee
and that $\pi>0$.
We shall refer to such invariant distribution $\pi$ as the \emph{centrality vector} of the social network.
The centrality vector characterizes the \emph{social power} of each agent in the group~\citep{French:1956,Friedkin:1986,ProTempo:2017-1} and has many other applications, see e.g.~\cite{Como.Fagnani:2015}.
The asymptotic behavior of the opinion dynamics \eqref{eq.degroot} can then be characterized as follows~\citep{Como.Fagnani:2016}.
\begin{proposition}\label{prop:limit}
Consider a social network with irreducible aperiodic influence matrix $P$ whose centrality vector is $\pi$.
Let  the opinion profile $x(t)$ evolve according to the French-DeGroot dynamics \eqref{eq.degroot} with self-confidence vector $z\in\mc Z$.
Then, a stochastic matrix  $H(z)$ exists such that
\begin{equation}\label{lim}\lim_{t\to+\infty}x(t)=H(z)x(0)\,,\end{equation}
for every initial configuration $x(0)$ in $\R^n$. Moreover,
\begin{enumerate}
\item[(i)] if $\mc S(z)=\emptyset$, then
\be\label{def:H(z)}H(z)=\frac1{\gamma(z)}\one\pi'(I-[z])^{-1}\,,\ee
where
\be\label{def:S(z)}\gamma(z)=\sum_{i=1}^n\frac{\pi_i}{1-z_i}\,;\ee
\item[(ii)] if $\mc S(z)\ne\emptyset$, then $H_{ij}(z)=0$ for every $i$ in $\mc V$,  $j$ in $\mc V\setminus\mc S(z)$, and $H_{ij}(z)=\delta_i^j$ for every $i$ and $j$ in $\mc S(z)$.
\end{enumerate}
\end{proposition}

\begin{remark}\label{rem:centrality} Notice that, when no agent is stubborn, $H(z)=\one\pi'$
has identical rows equal to $\pi'$, so that Proposition \ref{prop:limit} implies that every agent influences the others in a homogeneous way.
\end{remark}

\begin{remark}\label{rem:hitting} Upon interpreting the stochastic matrix $Q=(I-[z])P+[z]$ as the transition probability matrix of a discrete-time Markov chain, we have the following probabilistic interpretation. It is a known fact that, regardless of its initial state, such a Markov chain: (i) converges in probability to the invariant distribution $\gamma^{-1}(I-[z])^{-1}\pi$ if $\mc S(z)=\emptyset$; (ii) is absorbed in finite time in one of states from $\mc S(z)$ if $\mc S(z)\ne\emptyset$.
In the case (ii), $H_{ij}(z)$ is the probability of absorption in state $j$ when starting in state $i$.
\end{remark}

\subsection{Wisdom of agents and distributed estimation problem}

We assume that the initial opinions are the agents' guesses about some unknown state of the world, modeled as a scalar parameter $\theta$ in $\R$. More precisely, let
$$x_i(0)=\theta+\xi_i,\quad i\in\mc V\,,$$
 $\xi_i$, for $1\le i\le n$, are uncorrelated random noise variables with expected value $0$ and variance $\sigma_i^2>0$. The value $\sigma_i^{-2}$ may be interpreted as the ``wisdom'' of agent $i$: the ``wiser'' agent (i.e., the smaller $\sigma_i^2$),  the closer her initial opinion is to the actual state of the world $\theta$.

For brevity, we introduce the column vector of variances
$$
\sigma^2=(\sigma_1^2,\ldots,\sigma_n^2)\,.
$$

It then follows from Proposition \ref{prop:limit} that the asymptotic opinion of agent $i$ in $\mc V$ is given by
$$\lim_{t\to+\infty}x_i(t)=\theta+\zeta_i\,,$$
where
$$\zeta_i\dfb\sum_{j=1}^nH_{ij}(z)\xi_j\,,$$
represents the remaining noise in agent  $i$'s estimation of the state of the world $\theta$ at the end of the social interaction.
Such remaining noise $\zeta_i$ results in a convex combination of the noises $\xi_j$ in the initial agents' observations $x_j(0)$. As such, it has expected value $0$ and, since the $\xi_j$'s are assumed to be uncorrelated, its variance is given by\be\label{estimation-cost}\ups_i(z)=\E[\zeta_i^2]=\sum_{j=1}^nH_{ij}^2(z)\sigma_j^2\,.\ee

We shall interpret the above as the asymptotic estimation cost of agent $i$.
Given the influence matrix $P$ and the agents' initial variance values $\sigma_i^2$, such estimation cost depends only on the self-confidence vector $z$. Henceforth, we shall adopt a common abuse of notation writing
$$\ups_i(z)=\ups_i(z_i,z_{-i})$$
to emphasize the fact that the asymptotic estimation cost of agent $i$ depends both on her own self-confidence value $z_i$ and on  all other agents' self-confidence profile
$$z_{-i}=(z_j)_{j\ne i}\,.$$

Assume that the agents are rational entities selecting their self-confidence values in order to optimize the accuracy of their estimation of true state of the world $\theta$.

We shall use the following multi-objective optimization and  game-theoretic notions.
\begin{definition}
A self-confidence profile $z$ in $\mc Z$ is
\begin{enumerate}
\item[(i)] \emph{Pareto-optimal} if there exists no other self-confidence profile $\ov z$ in $\mc Z$ such that
$$\ups_i(\ov z)\le\ups_i(z)\,,\quad\forall i\in\mc V\,,$$ and there exists an agent $j$ in $\mc V$ such that $$\ups_j(\ov z)<\ups_j(z)\,;$$
\item[(ii)] a \emph{Nash equilibrium} if for every agent $i$ in $\mc V$
$$\ups_i(z_i,z_{-i})\le\ups_i(\ov z_i,z_{-i})\,,\qquad \forall\ov z_i\in[0,1]\,;$$
\item[(iii)] a \emph{strict Nash equilibrium} if for every agent $i$ in $\mc V$
$$\ups_i(z_i,z_{-i})<\ups_i(\ov z_i,z_{-i})\,,\qquad \forall\ov z_i\in[0,1]\setminus\{z_i\}\,.$$
\end{enumerate}
We shall denote by $\mc P\subseteq\mc Z$ the Pareto frontier, i.e., the set of Pareto optimal self-confidence profiles and by $\mc N$ and $\mc N^*$, respectively, the set of Nash equilibria and the set of strict Nash equilibria.
\end{definition}

\section{Main results}\label{sec:main}

This section presents main results of the paper. 

Our first main result characterizes the Pareto frontier and the set of strict Nash equilibria.
\begin{theorem}\label{theo:main}
Consider a social network with irreducible aperiodic influence matrix $P$ and centrality vector $\pi$.
Let the initial opinions of the agents be uncorrelated with variance $\sigma_i^2>0$ for $i$ in $\mc V$ and $\mc Z^*$ denote the set
\be\label{Z*}\mc Z^*=\{\one-\alpha[\pi]\sigma^2:\,0<\alpha\leq\alpha_*\},\;\; \alpha_*=\frac{1}{\max_i\{\pi_i\sigma_i^2\}}.\ee
Then, $\mc P=\mc N^*=\mc Z^*\,.$
\end{theorem}

Our second main result provides a characterization of the set of non-strict Nash equilibria.
In order to state it, recall the definition of restricted graph $\mc G[\mc S]$ from Section \ref{sec:notation}.

\begin{theorem}\label{theo:other-Nash}
Consider a social network with irreducible aperiodic influence matrix $P$ and centrality vector $\pi$.
Let the initial opinions of the agents be uncorrelated with variance $\sigma_i^2>0$ for $i$ in $\mc V$.
Then,
every possible non-strict Nash equilibrium $z$ in $\mc N\setminus\mc N^*$ is such that:
\begin{enumerate}
\item[(i)] $|\mc S(z)|\geq 2$;
\item[(ii)] $\sigma_i^2=\sigma_j^2$ for every $i$ and $j$ in $\mc S(z)$;
\item[(iii)] the graph $\mc G_P[\mc S(z)]$ is a directed ring.
\end{enumerate}
\end{theorem}

Theorem \ref{theo:other-Nash} has the following direct consequences.
\begin{corollary}\label{coro:no-nonstrct} If $\sigma_i^2\neq\sigma_j^2$ for $i\neq j$ in $\mc V$, then $\mc N=\mc Z^*$.
\end{corollary}
\begin{pf} By Theorem \ref{theo:other-Nash}, there are no non-strict Nash equilibria under the given assumptions: $\mc N=\mc N^*$. The claim now follows from Theorem~\ref{theo:main}.  \qed\end{pf} \medskip

 \begin{corollary}\label{coro:undirected} If $\mc G_P$ is undirected, then $|\mc S(z)|=2$ for every possible non-strict Nash equilibrium $z$ in $\mc N\setminus\mc N^*$.
\end{corollary}
\begin{pf}  Condition (iii) in Theorem \ref{theo:other-Nash} is never satisfied if $\mc G_P$ is undirected and $|\mc S(z)|\ge3$.\qed\end{pf}\medskip

While Theorem \ref{theo:other-Nash} and its corollaries provide necessary conditions for the existence of non-strict Nash equilibria, the following provides a sufficient condition.
\begin{proposition}\label{prop:main} Consider a social network with irreducible aperiodic influence matrix $P$ and uncorrelated agents' initial opinions with variances $\sigma_i^2$, $i\in\mc V$.
Then, every self-confidence vector $z\in\mc Z$ such that \be\label{cond}\sigma^2_j\le\sigma_i^2\,,\qquad  \forall j\in\mc S(z)\,,\ \forall i\in\mc V\,,\ee
and $|\mc S(z)|=2$ is a non-strict Nash equilibrium.
\end{proposition}

Notice that condition~\eqref{cond} entails that $\sigma^2_j=\sigma^2_k$, for every $j$ and $k$ in $\mc S(z)$.

\subsection{Numerical example}\label{subsec:exampl}

\begin{figure}[ht]
\centering
\includegraphics[width=0.75\columnwidth]{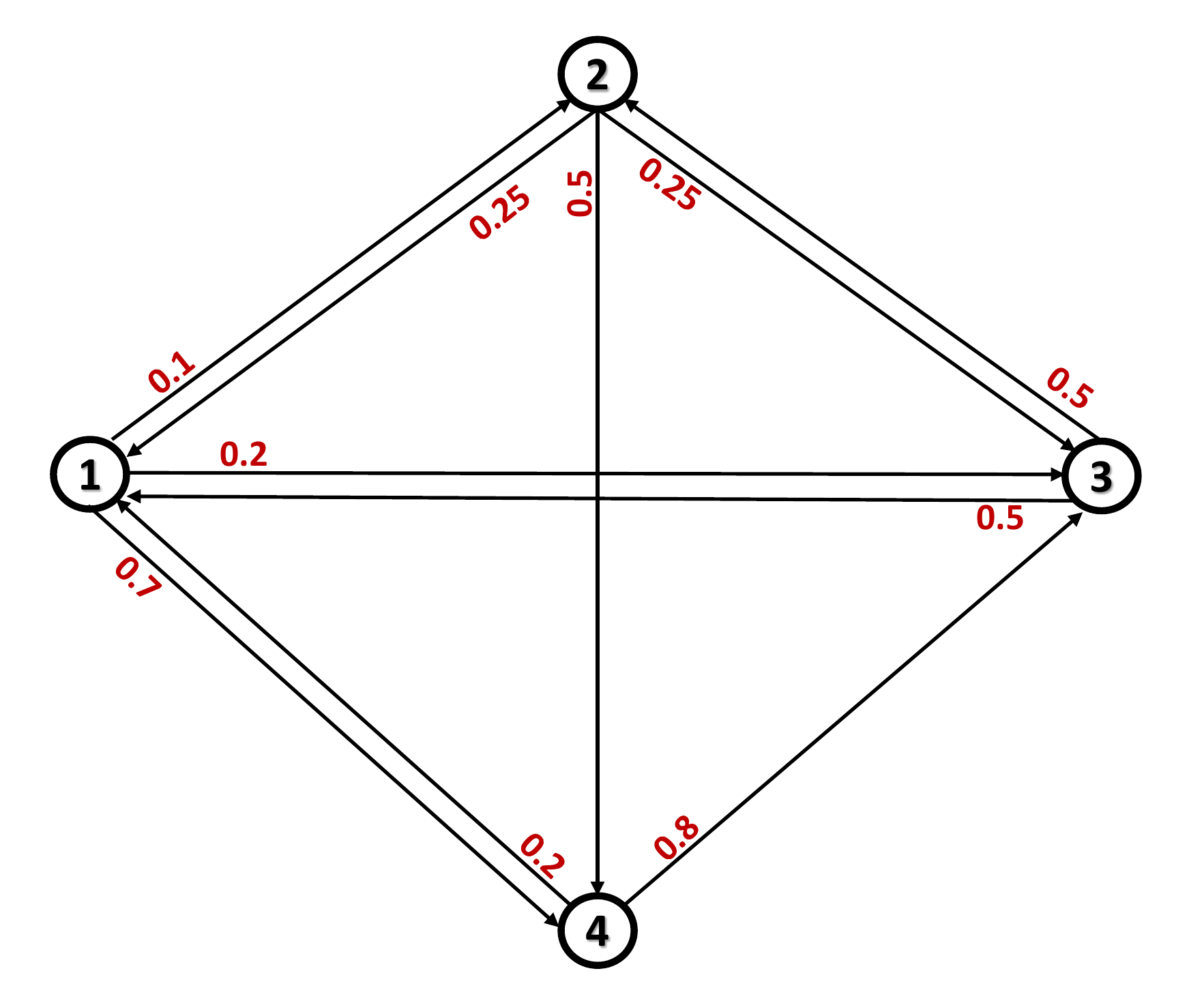}
\caption{Graph $\mc G_P$ of~\eqref{eq.p-example}}\label{fig.graph}
\end{figure}

Consider a social network with $n=4$ agents and influence matrix $P$ as follows
\begin{equation}\label{eq.p-example}
  P=\begin{bmatrix}
      0 & 0.1 & 0.2 & 0.7\\
      0.25 & 0 & 0.25 & 0.5\\
      0.5 & 0.5 & 0 & 0\\
      0.2 & 0 & 0.8 & 0
    \end{bmatrix}
\end{equation}
whose (weighted) graph is shown in Figure \ref{fig.graph}.
Solving equation~\eqref{pi}, the centrality vector $\pi$ is found to be
$$\pi=(0.2507,\; 0.1783,\; 0.3064,\; 0.2646)\,.$$
If the agents had no self-influence ($z=0$), agent $3$ would be most influential, whereas agent $2$ would have the least social power among all agents.

Let the vector of variances be
$$
\sigma^2=(0.32^2,\;0.35^2,\; 0.38^2,\; 0.29^2),
$$
that is, agent $4$ is the ``wisest'' and agent $3$ is the most ``unwise''.
In view of Theorem~\ref{theo:main} and Corollary~\ref{coro:no-nonstrct}, the set of Pareto optimal self-confidence profiles (and the set of all Nash equilibria) are given by~\eqref{Z*}, where
$$[\pi]\sigma^2=(0.0257,\;0.0218,\;0.0442,\;0.0223).$$
Notice that the largest self-confidence value in every Pareto optimal self-confidence value profile is  then agent $2$'s, who is neither most influential in terms of the eigenvector centrality nor the wisest. On the other hand,  the most ``powerful'' agent $3$ has the smallest self-confidence value due to its large variance.

\section{Proofs}\label{sec.proofs}
In this section, we develop our technical analysis leading to the proofs of our main results: Theorem \ref{theo:main}, Theorem \ref{theo:other-Nash}, and Proposition \ref{prop:main}.

\subsection{The best response correspondence and its properties}
We write $\mathcal B_i(z_{-i})\subseteq [0,1]$ to denote the best response (BR) set of an agent $i$ when the rest of agents have chosen self-confidence profile $z_{-i}$, i.e.,
\[
\mathcal B_i(z_{-i})=\argmin_{z_i\in[0,1]}\ups_i(z_i,z_{-i}).
\]

Our study starts with the case when all agents are playing $z_i<1$ or when there is just one agent playing $z_i=1$.

\begin{lemma}\label{lem.brm}
The BR set of agent $k$  in $\mc V$ satisfies
\beq\label{eq.brm}
\mc B_k(z_{-k})=\left\{\left[1-\frac{a_{-k}\pi_k\sigma_k^2}{b_{-k}}\right]_+\right\}\,,\quad\forall 0\le z_{-k}<\one\,
\eeq
where
\beq\label{eq.aux1}
a_{-k}
=\sum_{j\ne k}\frac{\pi_j}{1-z_j},\qquad
b_{-k}
=\sum_{j\ne k}\frac{\pi_j^2\sigma_j^2}{(1-z_j)^2}.
\eeq
\end{lemma}
\begin{pf}
For $0\le z<\one$, it is convenient to introduce new variables
\beq
y_i\dfb\frac{1}{1-z_i}\in [1,+\infty),\qquad i\in\mc V\,,
\eeq
and rewrite the cost of agent $k$ as
\[
\ups_k(z)=\ov\ups_k(y)=
\frac{(\pi_k^2\sigma_k^2y_k^2+b_{-k})}{(\pi_ky_k+a_{-k})^2}.
\]
Differentiating the function above with respect to $y_k$ yields
\be\label{eq.partial}
\ba{rcl}
\frac{\partial\ov\ups_k}{\partial y_k}(y)
&=&\frac{2\pi_k^2\sigma_k^2y_k}{(\pi_ky_k+a_{-k})^2}-
\frac{2(\pi_k^2\sigma_k^2y_k^2+b_{-k})\pi_k}{(\pi_ky_k+a_{-k})^3}
\\[10pt]
&=&\frac{2\pi_k(a_{-k}\pi_k\sigma_k^2y_k-b_{-k})}{(\pi_ky_k+a_{-k})^3}.
\ea
\ee
Hence, for $y_{-k}$ being fixed, $y_k\mapsto\ov\ups_k(y_k,y_{-k})$ achieves a strict minimum point at $y_k^*=\min\{1,b_{-k}/(a_{-k}\pi_i\sigma_k^2)\}$. This minimizer corresponds to
\[
z_k^*= 1-\frac{1}{y^*_k}=\left[1-\frac{a_{-k}\pi_i\sigma_k^2}{b_{-k}}\right]_+\,.
\]

Now, observe that, for  every fixed $0\le z_{-k}<\one$, Proposition \ref{prop:limit} implies that
$$\lim_{z_k\uparrow 1}H(z_k,z_{-k})=\one(\delta^k)^T=H(1,z_{-k})\,,$$
so that,  by \eqref{estimation-cost}, the function  $z_k\mapsto\ups_k(z_k,z_{-k})$ is continuous on $[0,1]$. This implies that $z_k^*$ is in fact the unique minimizer of $z_k\mapsto\ups(z_k,z_{-k})$ over closed interval $[0,1]$.\hfill\qed
\end{pf}

In the case where $z_{-k}$ contains some entries equal to $1$ (whose corresponding agents are stubborn), the analysis of the best respoense $\mc B_k(z_{-k})$ is more involved. We start with the following result.

\begin{lemma}\label{lem.brm1}
For an agent $k$ in $\mc V$ and $0\le z_{-k}\le\one$ such that $\mc S(z_{-k})\neq\emptyset$, the statements hold:
\begin{enumerate}
\item[(i)] $z_k\mapsto\ups_k(z_k,z_{-k})$ is constant on $[0,1)$; \\
\item[(ii)] $\mc B_k(z_{-k})\in\{[0,1), [0,1], \{1\}\}$\\
\item[(iii)] if $\sigma_k^2>\sigma_i^2$ for all $i$ in $\mc S(z_{-k})$, then $\mc B_k(z_{-k})=[0,1)$.
\end{enumerate}
\end{lemma}
\begin{pf}
Recall (c.f.~Remark \ref{rem:hitting}) that $H_{ki}(z)$ is the probability that the dissctrete-time random walk with transition probability matrix $(I-[z])P+[z]$ that starts in $k$ gets absorbed in $i$. This does not depend on the value $z_k\in [0,1)$ that only affects the exit time from node $k$. This proves (i), which trivially implies (ii).
Statement (iii) is a direct consequence of Proposition~\ref{prop:limit}: under the assumption of (iii),
for each $z_k\in[0,1)$ one has
$v_i(z_k,z_{-k})=\sum\nolimits_{j\in S(z_{-k})}H_{ij}^2(z)\sigma_j^2<\sigma_k^2=u(1,z_{-k})$.
\qed
\end{pf}

\textbf{Proof of Proposition \ref{prop:main}.}

As direct application of Lemma \ref{lem.brm1} we have the following proof of Proposition \ref{prop:main}.
It follows from \eqref{estimation-cost}, Proposition \ref{prop:limit} (ii), and \eqref{cond} that
$\sigma_j^2=\sigma_{\min}^2=\min\{\sigma_k^2:\,k\in\mc V\}\,\forall j\in\mc S(z)$, and therefore
$$\ups_i(z)=\sum_{j\in\mc S(z)}H^2_{ij}(z)\sigma_j^2\le\sigma_{\min}^2\le\sigma_i^2$$
for every $i$ in $\mc V$ with equality for $i$ in $\mc S(z)$. It then follows from Lemma \ref{lem.brm1} that
$\mc B_i(z_{-i})\supseteq[0,1)\ni z_i$ for every $i$ in $\mc V\setminus\mc S(z)$ and
$\mc B_j(z_{-j})=[0,1]$ for every $j$ in $\mc S(z)$, thus implying that $z$ is a non-strict Nash equilibrium.
\qed

\subsection{The structure of NE in the game}
The following simple result is proved in \cite{Peluffo:2019}.
\begin{lemma}\label{lemma:min}
For positive variance values $(\sigma_i^2)$, the convex  program
\be\min_{\substack{\mu\in\R_+^n:\\\one'\mu=1}}\sum_i\mu_i^2\sigma_i^2\,,\ee
has a unique solution
\be\mu^*_i=\frac{\sigma_i^{-2}}{\sum_j\sigma_j^{-2}}\,.\ee
\end{lemma}

We are now ready to prove our main resuls.

\textbf{Proofs of Theorem~\ref{theo:main} and Theorem~\ref{theo:other-Nash}}

It follows from Proposition \ref{prop:limit} and Lemma \ref{lemma:min} that, for every agent $i$ in $\mc V$,
\be\label{optimal}\ups_i(z)=\sum\nolimits_jH_{ij}^2(z)\sigma_j^2\ge\sum\nolimits_j(\mu_j^*)^2\sigma_j^2\,,\ee
with strict inequality unless $H_{ij}(z)=\mu_j^*$ for every $j$ in $\mc V$. This directy implies that the Pareto frontier is given by\footnote{Here, we some abuse of notation, $\sigma^2$ is used to denote vector $(\sigma_i)^2$.}
\be\label{eq.p=z}
\ba{rcl}\mc P&=&\{z\in\mc Z:\,H(z)=\one(\mu^*)^{\top}\}\\
&=&\{z\in[0,1)^n:\,(I-[z])^{-1}\pi=\gamma(z)\mu^*\}\\
&=&\{z\in[0,1)^n:\,[\mu^*]^{-1}\pi=\gamma(z)(\one-z)\}\\
&=&\{z=\one-\alpha[\pi]\sigma^2:\,0<\alpha<\alpha_*\}\\
&=&\mc Z^*\,.
\ea
\ee
On the other hand, it directly follows from \eqref{optimal} that every $z$ in $\mc Z^*$ is a strict Nash equilibrium.

We now prove that there are no Nash equilibria in the set $[0,1)^{n}\setminus\mc Z^*$.
Indeed, assume that $0\le\bar z<\one$ is a Nash equilibrium. By a change of variables, we can replace $P$ with $(I-[\ov z])P+[\ov z]$, $\pi$ with $\ov\pi=\gamma(\ov z)^{-1}(I-[\ov z])^{-1}\pi$ and $z$ with $z+\bar z$.
We now introduce auxiliary variables as in the proof of Lemma~\ref{lem.brm}: $y_i=(1-z_i)^{-1}$ for $i\in\mc V$. We set
\beq
A
=\sum_{j}\frac{\ov\pi_j}{1-z_j},\qquad
B
=\sum_{j}\frac{\ov\pi_j^2\sigma_j^2}{(1-z_j)^2}\eeq
and we obtain, using derivation \eqref{eq.partial} and expressions \eqref{eq.aux1},
$$\frac{\partial\ups_k}{\partial y_k}(y)=\frac{2\pi_k(A\ov\pi_k\sigma_k^2y_k-B)}{(\ov\pi_ky_k+a_{-k})^3}
$$
We recall that, by our assumptions, $z=0$ (and thus $y=\one$) must correspond to a NE. Consequently, for $y=\one$ all partial derivatives must be non positive. We notice that the numerators, for $y=\one$, satisfy the following property:
$$\sum\limits_{k\in\mc V}2\ov\pi_k(B-A\ov\pi_k\sigma_k^2)=B-A\sum\limits_k\ov\pi_k\sigma_k^2=0$$
where last equality follows from the special form of $A$ and $B$ for $z=0$. Consequently, we must have that
$$B-A\ov\pi_k\sigma_k^2=0,\quad\forall k\in\mc V$$
This says that $\ov\pi$ is proportional to vector $(\sigma_i^{-2})$, and thus $\ov\pi=\mu^*$, which is equivalent to $\ov z\in\mc Z^*$ due to~\eqref{eq.p=z}.

We are left with the analysis of possible NE $z$ for which $S(z)\neq\emptyset$.
First, notice that for every self-confidence profile $z$ in $\mc Z$ such that $S(z)=\{k\}$, Lemma \ref{lem.brm} implies $1\notin\mc B_k(z_{-k})$. This shows that $|\mc S(z)|\ne 1$ for every Nash equilibrium $z$, thus implying the necessity of condition (i) in Theorem \ref{theo:other-Nash}.

On the other hand, let $z$ in $\mc Z$ be a self-confidence profile such that  $|S(z)|\geq 2$ and let
$S^*(z)=\argmax_{i\in\mc S(z)} \sigma_i^2$. If $S^*(z)\neq S(z)$, 
then in $\mc G_P$, there must exist agents $k$ in $\mc S^*(z)$ and $j$ in $\mc S\setminus\mc S^*(z)$ and
a sequence
$i_1,\ldots,i_{l-1},i_l$ such that
$$i_0=k\,,\quad i_l=j\,,\quad  i_s\not\in S(z)\;\forall 1\le s<l\,,$$
and
$$\prod_{1\le  s< l}(1-z_{i_s})P_{i_{s}i_{s+1}}>0\,.$$
It then follows from Proposition \ref{prop:limit} that, for every $\tilde z$ such that $0\le\tilde z_k<1$ and $\tilde z_{-k}=z_{-k}$,
$$H_{kj}(\tilde z)\ge\prod_{0\le  s< l}(1-z_{i_s})P_{i_{s}i_{s+1}}>0\,,$$
so that
$$\ups_k(\tilde z)<\ups_k(z)=\sigma_k^{-2}\,,$$
thus proving that $z$ cannot be a Nash equilibrium.
This implies that $S^*(z)= S(z)$, thus proving condition (ii) of Theorem \ref{theo:other-Nash}.

Finally, consider a non-strict Nash equilibrium $z$. By points (i) and (ii), $|\mc S(z)|\ge2$ and $\sigma_j^2=\sigma_*^2$ for every $j$ in $\mc S(z)$. 
Now, recall the definition of the restricted graph $\mc G_P[\mc S(z)]=(\mc S(z),\mc E[\mc S(z)])$ and first notice that, since $\mc G_P$ is strongly connected, so is $\mc G_P[\mc S(z)]$. 
For $i$ in $\mc S(z)$, let  $\mc J_i=\{j\in\mc S(z):\,(i,j)\in\mc E[\mc S(z)]\}$ be the set of out-neighbors of node $i$ in $\mc G[\mc S]$ and observe that $H_{ij}(\tilde z_i,z_{-i})>0$ for every $\tilde z_i$ in $[0,1)$  if and only if $j\in\mc J_i$. Then, \eqref{estimation-cost} implies that
$$\ups_i(\tilde z_i,z_{-i})=\sum_{j\in\mc J_i}H_{ij}^2(z)\sigma_j^2=\sum_{j\in\mc J_i}H_{ij}^2(z)\sigma_*^2\le\sigma^2_*=\ups_i(z)\,,$$
with equality if and only if $|\mc J_i|=1$. This proves that $\mc G_P[\mc S(z)]$ is $1$-regular. Since it is also strongly connected, we conclude  that it is a directed ring, 
thus completing the proof of the result.
\hfill\qed

\section{Myopic Self-Confidence Adaptation Dynamics}\label{sec:dynam}



In this section, we introduce a family of continuous-time self-confidence adaptation dynamics and present some numerical simulations. The outcome of the presented simulations and other non-reported ones lead us to conjecture that such self-confidence adaptation dynamics globally converge to one of the equivalent strict Nash equilibria of the game. We shall present an analytical study addressing this conjecture in a subsequent work. We also believe that asynchronous discrete-time versions of such self-confidence adaptation dynamics can be studied and may display analogous asymptotic behavior.

The considered dynamics are obtained by equating the change rate of logarithmic derivative of the self-confidence of an agent, i.e., $\frac{d}{dt}{(\log z_i)}=\dot z_i/z_i$ to the gradient of her utility, as follows,
\begin{equation}\label{eq.dynam}
  \dot z_i=-z_i\frac{\partial v_i(z_i,z_{-i})}{\partial z_i},\quad i=1,\ldots,n.
\end{equation}
The main consequence for the presence of the term $z_i$ in the right hand side of the ODE is to avoid vector fields to point out of the self-confidence profile space $\mc Z$ on boundary points where some of the actions are equal to $0$. In fact, it can be shown that the set of self-confidence profiles $\{0\leq z<\one\}$ is forward-invariant for \eqref{eq.dynam}.
Furthermore, numerical simulations show that, in fact, solutions starting in the open hypercube $\{0<z<\one\}$ in fact
do not converge to the boundary but rather converge to one of the strict Nash equilibria~\eqref{Z*},
which in turn can be shown to be equilibrium points of the system~\eqref{eq.dynam}.

We illustrate the behavior of the self-confidence adaptation dynamics~\eqref{eq.dynam} for the network of $n=4$ agents introduced in Section.~\ref{subsec:exampl} and different initial conditions.
Notice that $\alpha_*$ from Theorem~\ref{theo:main} is
$1/0.0442\approx 22.624$.

\textbf{Case 1.} The initial self-confidence values are large: $z_i(0)=0.994$ for all $i$. In this situation, the steady vector $z^*$ is very close to $z(0)$ and corresponds to the element of $\mathcal{Z}^*$ with $\alpha\approx 0.37$.

\textbf{Case 2.} The initial self-confidence values are relatively small: $z_i(0)=0.1$,  for all $i$. In this situation, agents $1$, $2$, and $4$ increase their self-confidence values substantially,
whereas the self-confidence value of agent $3$ becomes very small yet positive (about $2\cdot 10^{-4}$). The final value corresponds to the element of $\mathcal{Z}^*$ with $\alpha\approx 22.59$ (almost maximal value).

\textbf{Case 3.} The self-confidence value of agent $1$ is very small, and the weights of agents $2$ and $4$
are almost maximal: $z_1(0)=0.01, z_2(0)=z_4(0)=0.99$, $z_3(0)=0.5$. In this situation, all agents change their self-confidence values.
The self-confidence values of agents $2$ and $4$ remain almost equal and both decrease, the self-confidence value of agent $1$ increases dramatically.
Remarkably, the self-confidence value of agent $3$ is a non-monotone function of time whose steady value is less than its initial value.
The steady vector $z^*$ corresponds to $\alpha\approx 8.7$.

\textbf{Case 4.} In this scenario, the self-confidence values of four agents were sampled randomly from the uniform distribution. In this situation, two self-confidence values are non-monotone functions of time. The steady value corresponds to
$\alpha\approx 21.06$.

\begin{figure}[h!]
  \centering
  \begin{subfigure}[t]{\columnwidth}
  \centering
  \includegraphics[width=0.76\columnwidth]{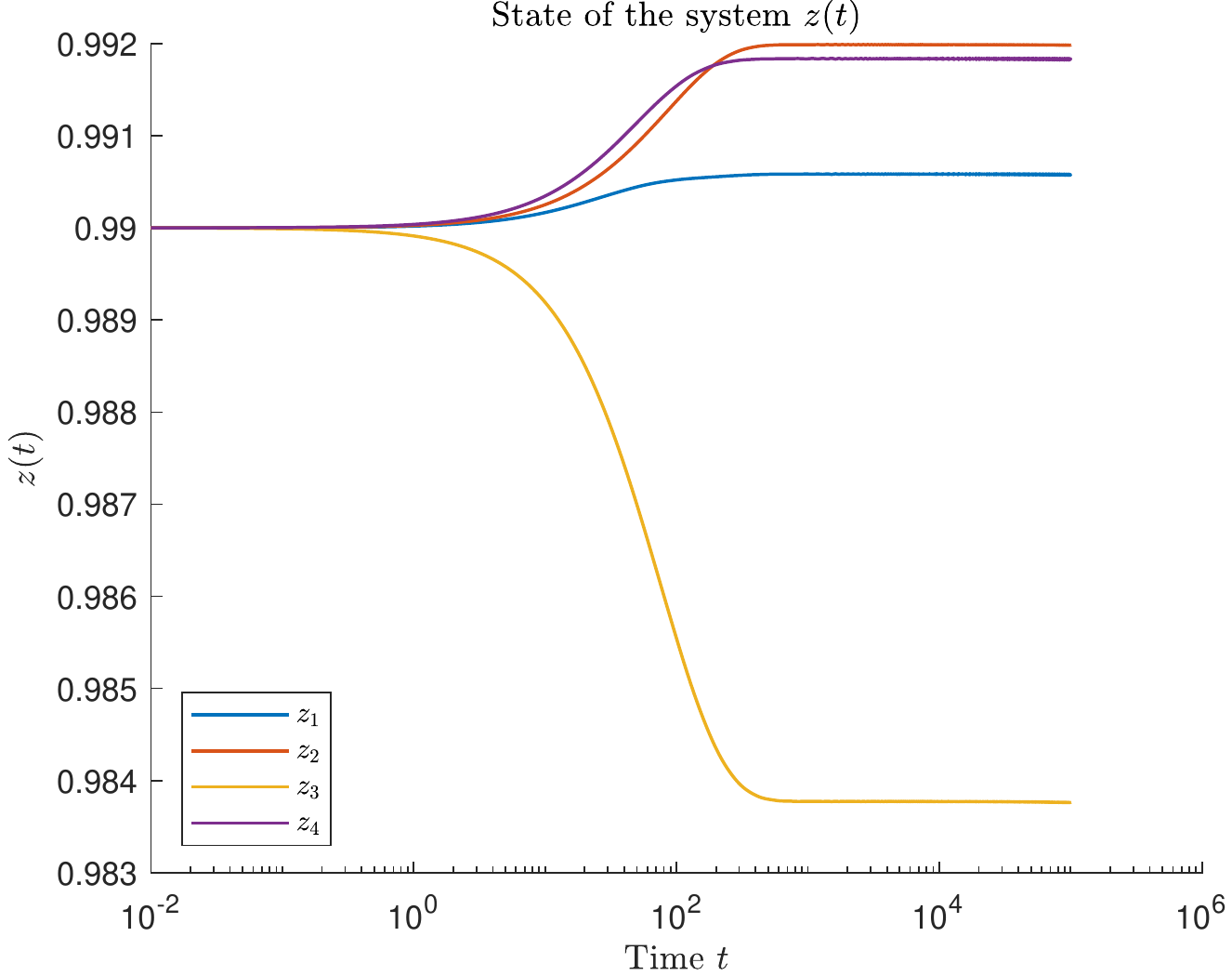}
  \end{subfigure}\\[2mm]
  \begin{subfigure}[t]{\columnwidth}
  \centering
  \includegraphics[width=0.76\columnwidth]{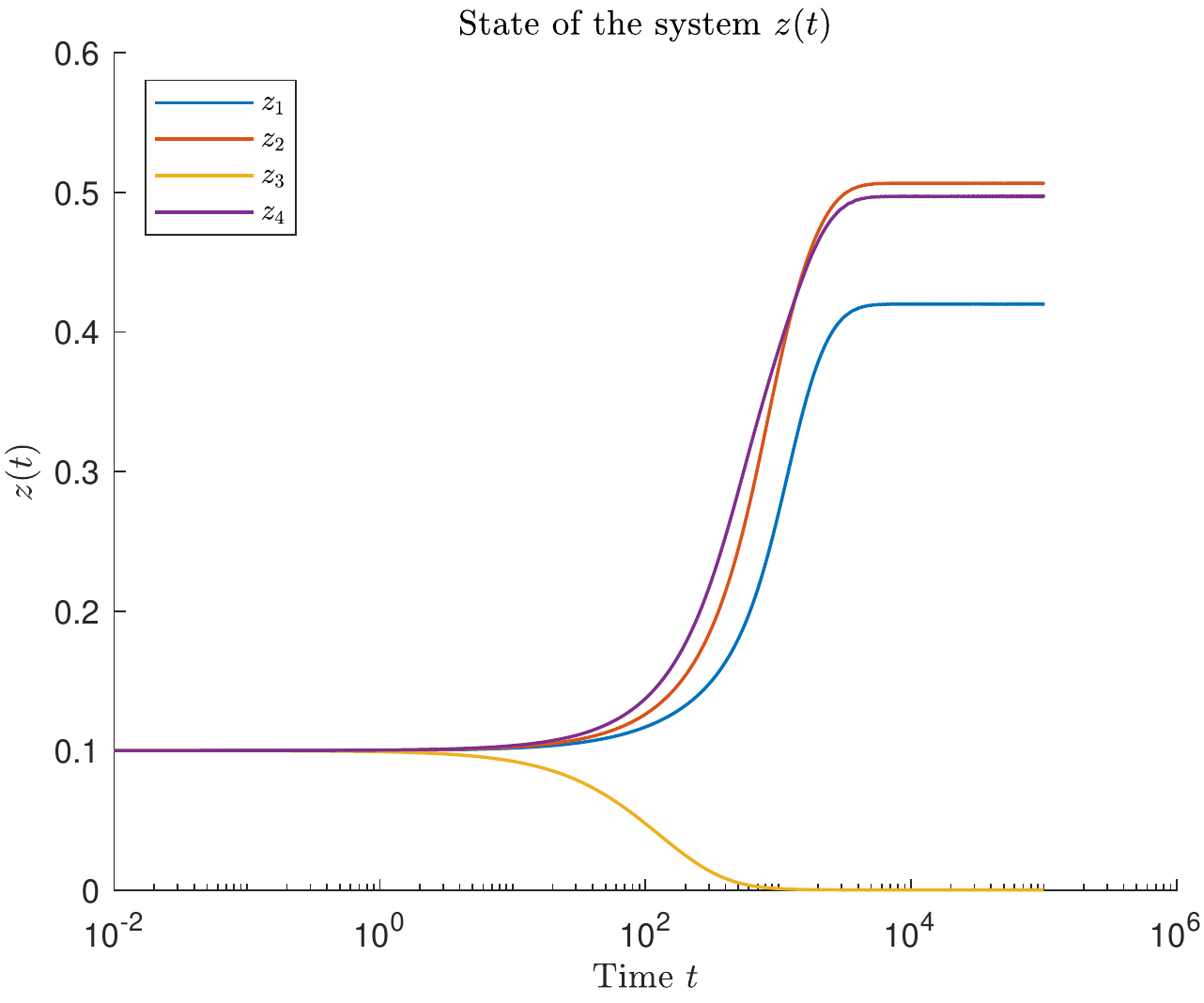}
  \end{subfigure}\\[2mm]
  \begin{subfigure}[t]{\columnwidth}
  \centering
  \includegraphics[width=0.76\columnwidth]{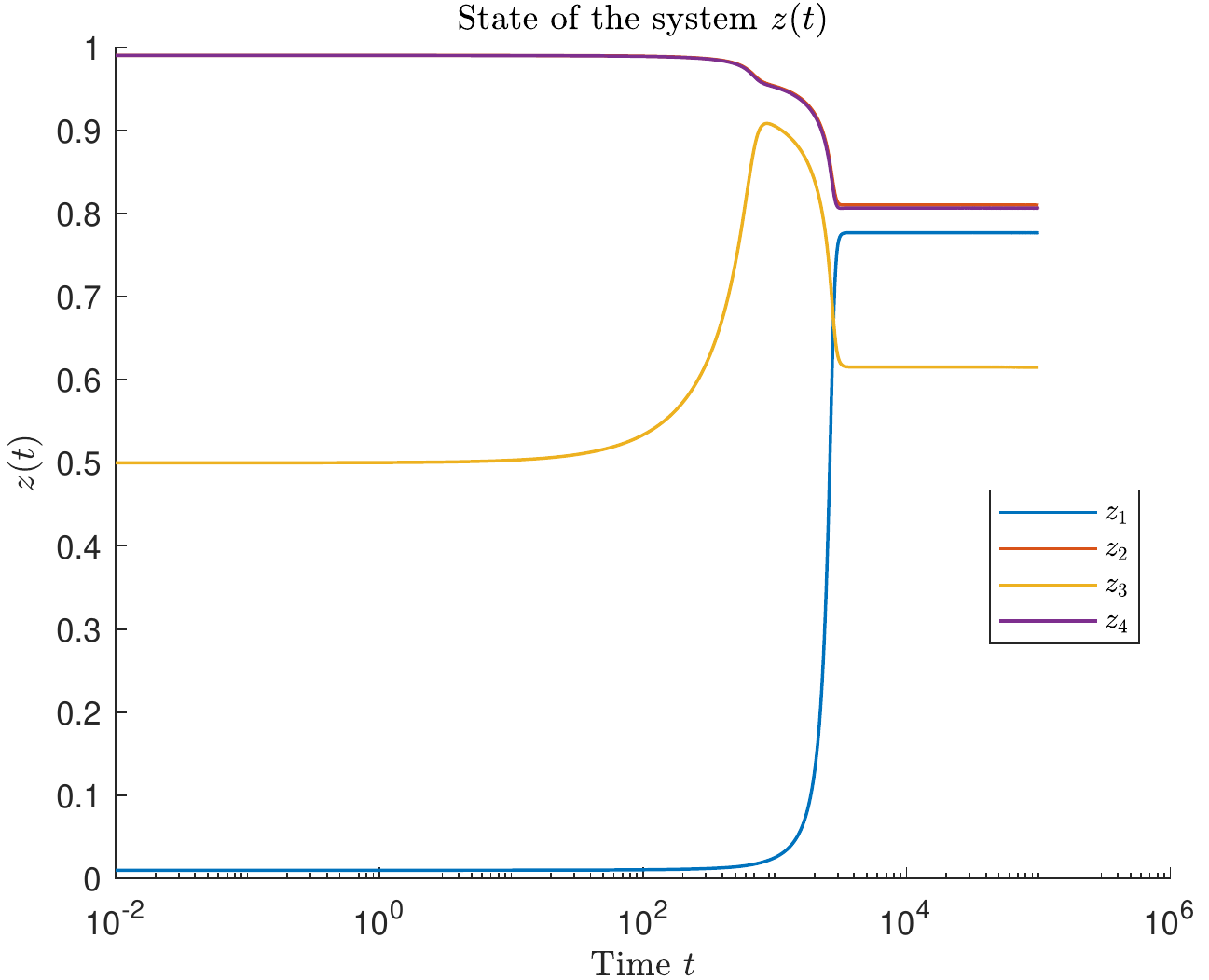}
  \end{subfigure}\\[2mm]
  \begin{subfigure}[t]{\columnwidth}
  \centering
  \includegraphics[width=0.76\columnwidth]{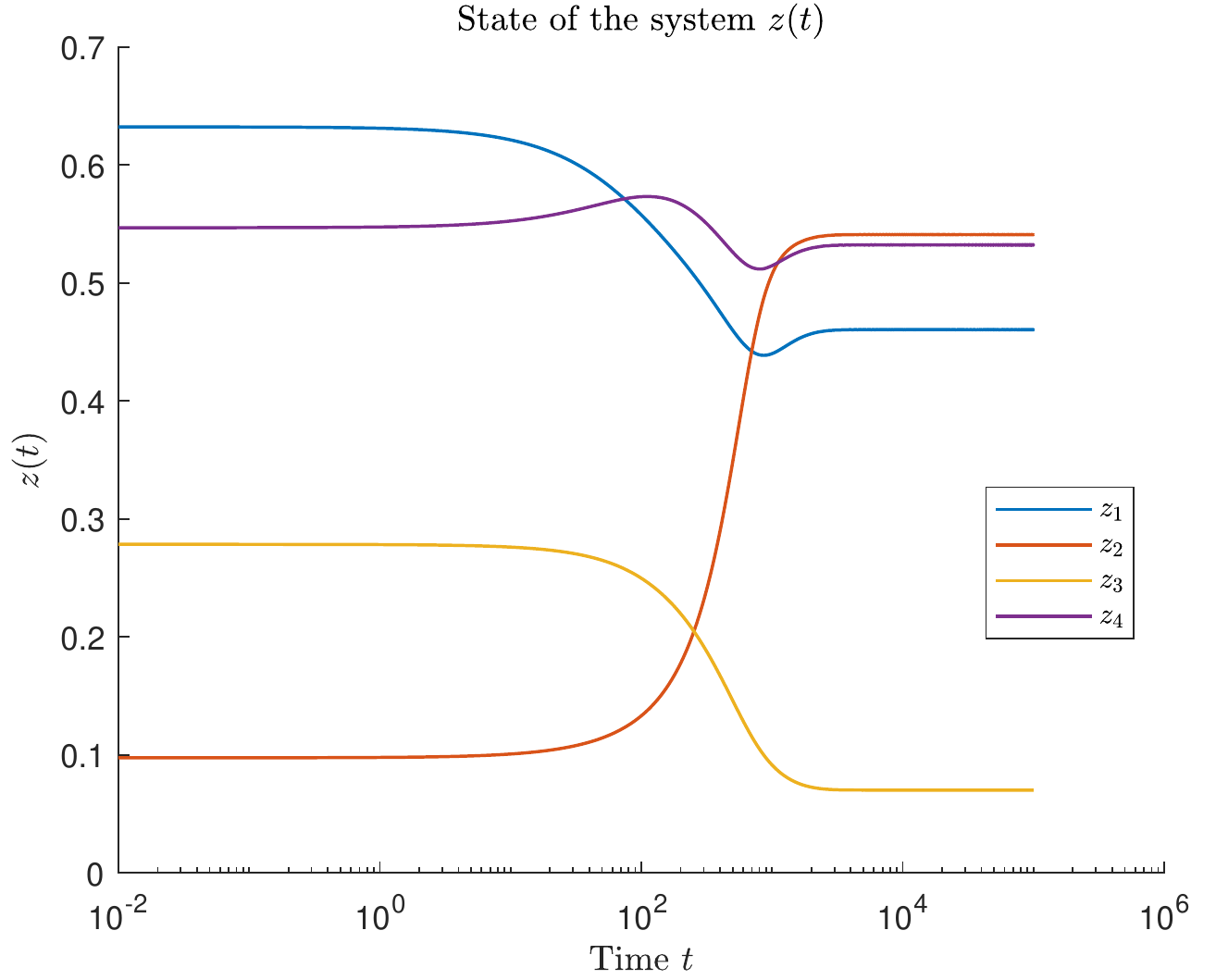}
  \end{subfigure}
  \caption{Solutions of the continuous-time synchronous myopic
self-confidence adaptation dynamics~\eqref{eq.dynam} for different initial self-confidence profiles $z(0)$ in $\mc Z$.}\label{fig.dynam}
\end{figure}

The results of simulations are shown in Figure \ref{fig.dynam}. One may notice that the convergence of continuous-time synchronous myopic self-confidence adaptation dynamics~\eqref{eq.dynam} appears to be rather slow.

\section{Future works}

We believe that there are two main directions for further study worth being pursued. On the one hand, the problem can be generalized, e.g., by allowing the agents to choose not just their self-confidence value, but also to redistribute the relative influence weights assigned to their neighbors in the network, or by assuming that the agents' objective is not merely minimizing their asymptotic estimate's variance, but also, e.g., increasing their centrality in the resulting network (c.f., \cite{Castaldo.ea:2020}). On the other hand, while the proposed analysis is purely static, it would be worth investigating the behavior of game-theoretic learning dynamics for the self-confidence evolution, such as best-response dynamics. One possible class of such dynamical systems is presented in Section~\ref{sec:dynam}.

Notice, however, that~\eqref{eq.dynam} is not a distributed algorithm computing Nash equilibria, because the gradients of cost functions depend on vectors $\pi$ and $\sigma^2$. This is consistent with the fact that we are considering a game whereby the utility functions of the players tend to depend on the actions of all other players and not just on those of their neighbors in the network.
Finding an efficient distributed algorithm is another topic of ongoing research.
\bibliography{mine,social}
\end{document}